\documentclass[letterpaper, 10 pt, conference]{ieeeconf}  
\usepackage{cite}
\usepackage{amsmath,amssymb,amsfonts}
\usepackage{algorithmic}
\usepackage{graphicx}
\usepackage{textcomp}
\usepackage{graphicx}
\usepackage[T1]{fontenc}
\usepackage{multicol}
\usepackage[ruled,linesnumbered]{algorithm2e}
\usepackage{textcomp}
\usepackage{xcolor}
\usepackage{soul}

\IEEEoverridecommandlockouts                              

\title{\LARGE \bf
The Novel Adaptive Fractional Order Gradient Decent Algorithms Design Via Robust Control
}

\author{Jiaxu Liu$^{1}$ , Song Chen$^{2}$, Shengze Cai$^{3}$ and Chao Xu$^{4}$
\thanks{*This work was supported by the National Key Research and Development Program of China under Grant 2019YFB1705800, the National Natural Science Foundation of China number 61973270, the Zhejiang Provincial Natural Science Foundation of China number LY21F030003, and the Zhejiang Provincial Natural Science Foundation of China number LY19A010024. (\textit{Corresponding author: Chao Xu}.)}
\thanks{$^{1}$J. Liu is with the School of Mathematical Sciences, Zhejiang University, Hangzhou, Zhejiang 310027, China (e-mail: jiaxuliu@zju.edu.cn).}%
\thanks{$^{2}$S. Chen is with the School of Mathematical Sciences, Zhejiang University, Hangzhou, Zhejiang 310027, China (e-mail: math\_cs@zju.edu.cn).}%
\thanks{$^{3}$S. Cai is with the State Key Laboratory of Industrial Control Technology, Institute of Cyber-Systems and Control, Zhejiang University, Hangzhou, Zhejiang 310027, China (e-mail: shengze\_cai@zju.edu.cn).}%
\thanks{$^{4}$C. Xu is with the State Key Laboratory of Industrial Control Technology, Institute of Cyber-Systems and Control, Zhejiang University, Hangzhou, Zhejiang 310027, China, and also with Huzhou Institute of Zhejiang University, Huzhou, Zhejiang 313000, China.  (e-mail: cxu@zju.edu.cn).}%
}

\begin{document}

\maketitle
\thispagestyle{empty}
\pagestyle{empty}

\begin{abstract}

The vanilla fractional order gradient descent may oscillatorily converge to a region around the global minimum instead of converging to the exact minimum point, or even diverge, in the case where the objective function is strongly convex. To address this problem, a novel adaptive fractional order gradient descent (AFOGD) method and a novel adaptive fractional order accelerated gradient descent (AFOAGD) method are proposed in this paper. Inspired by the quadratic constraints and Lyapunov stability analysis from robust control theory, we establish a linear matrix inequality to analyse the convergence of our proposed algorithms. We prove that the proposed algorithms can achieve R-linear convergence when the objective function is $\textbf{L-}$smooth and $\textbf{m-}$strongly-convex. Several numerical simulations are demonstrated to verify the effectiveness and superiority of our proposed algorithms.

\end{abstract}

\section{INTRODUCTION}

Gradient descent (GD) is a classical and effective method for solving unconstrained  optimization problems. Because of its simple structure and easy implementation, it has been widely used in different fields, such as control theory\cite{b1} and deep learning\cite{b2}. To improve the performance of the optimization algorithm, a number of  variants have been proposed based on the traditional GD algorithm, including momentum gradient method\cite{b5}, Nesterov's accelerated method\cite{b6}, Adam\cite{b7} and so on. However, a zigzag is commonly observed in the iterative process when searching the optimum of some typical functions, which affects the performance of the optimization problem. 

On the other hand, the fractional calculus, which has been reported to have some good properties (e.g., long-term memory\cite{b8}), has become a promising tool for improving the performance of gradient descent. In recent years, there are a few researches focusing on the fractional order gradient descent, especially on the convergence of the algorithms. 
Pu et al.\cite{b9} proposed a fractional order gradient method and pointed out that the convergence point of fractional order descent method is not the minimum point in the traditional sense due to the non-local property of the fractional order differential. In \cite{b10}, a fractional order gradient descent method for solving the  convex functions was proposed,  but this method cannot guarantee the convergence to the minimum point for strongly convex functions. In addition, Wei et al. \cite{b11} truncated the higher order terms of Taylor series to design the iteration formula. Following~\cite{b11}, Liu et al.\cite{b12} proposed a quasi fractional order gradient descent method. However, there existed the curse of dimensionality for solving \textbf{Armijo criterion} for these approaches. In this paper, we are inspired by the adaptive GD algorithms and we aim to combine them with the fractional order gradient. An adaptive fractional order gradient descent (AFOGD) method and an adaptive fractional order accelerated gradient descent (AFOAGD) method are proposed accordingly. Compared to the existing methods, AFOGD and AFOAGD not only retain the advantage of long-term  memory from fractional calculation, but also ensure the convergence by adaptively adjusting the step size. 

In order to analyze the convergence behavior of the aforementioned GD methods, we follow the researches where the optimization problem is analyzed based on control theory. 
A. Bhaya et al.\cite{b13} pointed out that the iterative solution algorithm of linear and nonlinear equations can be described as a feedback control system. They also deduced the continuous and discrete forms of Newton Raftson and conjugate gradient descent method by choosing an appropriate Lyapunov function. More recently,  Lessard et al.\cite{b17} developed a positive semi-definite programming framework based on quadratic Lyapunov functions, leading to sufficient conditions for exponential  convergence of an algorithm when the objective function is strongly convex. In addition, Su et al.\cite{b18} demonstrated that the Nesterov's accelerated gradient method can be expressed as an AVD equation\cite{19}. A similar idea was also proposed in \cite{b19}, where it was reported that both momentum gradient method and Nesterov's gradient method can be regarded as a PI-controlled closed-loop feedback system. Following these literature, we introduce the ideas from robust control theory in this paper to analyze the convergence of the proposed GD methods.

The main contributions in this paper can be summarized as follows. 
(1) We combine adaptive algorithms with fractional order gradient, resulting in two algorithms called AFOGD and AFOAGD.
(2) Based on robust control theory, we establish a linear matrix inequality (LMI) to analyse the convergence. Subsequently, we  prove that the AFOGD and AFOAGD  can achieve $R-$linear convergence rate if the objective function is  $L-$smooth and $m$-strongly-convex. 
(3) The numerical simulations verify the effectiveness and superiority of the AFOGD and AFOAGD, supporting our theoretical findings in this paper.

The rest of this article is organized as follows. In Section 2, the basic definition of smooth and strongly-convex functions, the vanilla fractional order gradient descent, and the relation between optimization and dynamical system are presented. The AFOGD and AFOAGD are introduced in Section 3 and the convergence rate of them are also analysed. Three numerical examples are given in Section 4 to illustrate the effectiveness and superiority of the proposed algorithms. Eventually, Section 5 concludes the paper.

\section{Preliminaries}

\subsection{Smooth and strongly convex function}
\textbf{Definition $1$ (smoothness)}: A differentiable function $f: \mathbb{R}^n \rightarrow \mathbb{R}$ is $L$-smooth if for all $x, y \in \mathbb{R}^n$
\begin{equation}
    \|\nabla f(x)-\nabla f(y)\|_2 \leq L\|x-y\|_2.
\end{equation}
It is equivalent to that for all $x, y \in \mathbb{R}^n$
\begin{equation}
    f(y) \leq f(x)+\nabla f(x)^{\top}(y-x)+\frac{L}{2}\|y-x\|_2^2.
\end{equation}

\textbf{Definition $2$ (strong convexity)}: A differentiable function $f: \mathbb{R}^n \rightarrow \mathbb{R}$ is $m$-strongly convex  if for all $x, y \in \mathbb{R}^n$
\begin{equation}
    \quad m\|x-y\|_2^2 \leq(x-y)^{\top}(\nabla f(x)-\nabla f(y)). 
\end{equation}
 It also implies that for all $x, y \in \mathbb{R}^n$
\begin{equation}
    f(x)+\nabla f(x)^{\top}(y-x)+\frac{m}{2}\|y-x\|_2^2 \leq f(y).  
\end{equation}
To combine strong convexity and Lipschitz continuity in a single inequality, we note that $\nabla f$ also satisfies \cite{b20} 
\begin{equation}
\begin{aligned}
    & \frac{m L}{m+L}\|y-x\|_2^2+\frac{1}{m+L}\|\nabla f(y)-\nabla f(x)\|_2^2 \\
    & \leq(\nabla f(y)-\nabla f(x))^{\top}(y-x).
\end{aligned}
\end{equation}
The above inequality can be represented by the following quadratic constraint:
\begin{equation}
    \begin{aligned}
\left[\begin{array}{c}
x-y \\
\nabla f(x)-\nabla f(y)
\end{array}\right]^{\top} Q_f\left[\begin{array}{c}
x-y \\
\nabla f(x)-\nabla f(y)
\end{array}\right] \geq 0,
\end{aligned}
\end{equation}
 where
 \begin{equation}
     \begin{aligned}
     Q_f=\left[\begin{array}{cc}
\frac{-m L}{m+L} I_n & \frac{1}{2} I_n \\
\frac{1}{2} I_n & \frac{-1}{m+L} I_n
\end{array}\right].
     \end{aligned}
 \end{equation}
 
\subsection{Fractional order gradient }

\textbf{Definition 3}: The definition of Caputo fractional-order
derivative for function $f(x)$ is given by \cite{b21}:
\begin{equation}
    { }_{x_0}^{C} \mathcal{D}_x^\mu f(x)=\frac{1}{\Gamma(t-\mu)} \int_{x_0}^x \frac{f^{(t)}(\tau)}{(x-\tau)^{\mu-t+1}} d \tau ,
\end{equation}
where $\mu \in (t-1, t)$ , $t\in N^+$ and $\Gamma(\cdot)$ is Gamma function.

The formula $(8)$ can be rewritten in Taylor series as follow \cite{b22}:
\begin{equation}
    { }_{x_0}^C \mathcal{D}_x^\mu f(x)=\sum_{i=t}^{\infty} \frac{f^{(i)}(x_0)}{\Gamma(i+1-\mu)}(x-x_0)^{i-\mu} .
\end{equation}
Similar to the general GD method, the fractional order gradient descent can be written as:
\begin{equation}
    x_{k+2}=x_{k+1}-\alpha_0 \cdot{ }_{x_{k}}^C \mathcal{D}_x^\mu f(x) .
\end{equation}
By substituting (9) into (10), we have
\begin{equation}
    x_{k+2}=x_{k+1}-\alpha_0 \sum_{i=0}^{\infty} \frac{f^{(i+1)}\left(x_{k}\right)}{\Gamma(i+2-\mu)}\left(x_{k+1}-x_{k}\right)^{i+1-\mu} .
\end{equation}
 However, it is almost impossible to compute the fractional order differential with the infinite series. Hence, we only reserve the first term in equation $(11)$, namely: 
 \begin{equation}
     { }_{x_{k}}^c \mathcal{D}_x^\mu f(x) \approx \frac{f^{\prime}\left(x_{k}\right)}{\Gamma(2-\mu)}\left(x-x_{k}\right)^{1-\mu} ,
 \end{equation}
where $\Gamma(2-\mu)$ becomes a constant. Then, the iteration formula (11) can be approximated by:
 \begin{equation}
      x_{k+2}=x_{k+1}-\alpha \cdot f^{\prime}\left(x_{k}\right)\left(x_{k+1}-x_{k}\right)^{1-\mu} ,
 \end{equation}
 with $0<\mu<1$ and $\alpha=\frac{\alpha_0}{\Gamma(2-\mu)}>0$ .
 The fractional order GD method can be directly extended to the case of $1<\mu<2$. We note that $\left(x_{k+1}-x_{k}\right)^{1-\mu}$ exists complex operations or is possibly infeasible. In order to solve the above problems, we modify equation $(13)$ to:
 \begin{equation}
     x_{k+2}=x_{k+1}-\alpha f^{\prime}\left(x_{k+1}\right)(||x_{k+1}-x_{k}||_2+\delta)^{1-\mu}
 \end{equation}
  which is called fractional order gradient descent (FOGD), where $0<\mu <2$ and $\delta>0$ is a constant\cite{b23}.
 
However, when $f(x) $ is a $m-$strongly convex function and $x_*$ is a unique fixed point of $f(x)$, if the iterative optimization calculation is performed using FOGD with $1<\mu<2$, the resulting sequence $ \left\{ x_k \right\}$ may diverge or converge to a region containing $x_*$ instead of converging to the true fixed point $x_*$. Below is an example. 

\textbf{Example 1.}  Let us consider a simple  quadratic function $f(x_1,x_2)=2x_1^2+3x_2^2+3$. It is easy to see that the fixed point of $f(x_1,x_2)$ is $x_*=[0,0]^{\top}$ and the minimum is $3$. The initial states are $x_0=[0.1,0.1]^{\top}$, $x_1=[1,1]^{\top}$. By selecting $\delta=10^{-4}$, $\alpha=0.2$ and $\mu=1.7$, the result is shown in Fig.1, where it can be seen that the FOGD method does not guarantee $f(x_1,x_2)$ converge to a minimum. This also motivates our research and raises the question: how to modify FOGD so that it can find the minimum for strongly convex functions?
 \begin{figure}
    \centering
    \includegraphics[width=\linewidth]{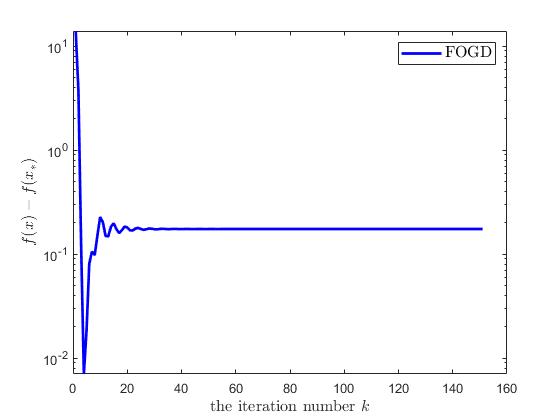}
    \caption{An example showing that the fractional order gradient descent method cannot find the minimum point for a quadratic function:  $f(x_1,x_2)=2x_1^2+3x_2^2+3$.}
    \label{fig:my_label}
\end{figure}

\subsection{Optimization and Control}
For a first-order optimization algorithm, it can be formulated as a dynamical system\cite{b24}:
\begin{equation}
    \begin{aligned}
\xi_{k+1} &=A \xi_k+B u_k \\
y_k &=C \xi_k \\
u_k &=\phi\left(y_k\right) \\
x_k &=E \xi_k
\end{aligned}
\end{equation}
where at each iteration (indexed by $k$), $u_k \in \mathbb{R}^n$ is the input $(n \leq d)$; $\xi_k \in \mathbb{R}^d$ and $y_k \in \mathbb{R}^n$ are the state and the feedback output, respectively; $\phi:R^n\rightarrow R^n $ represents the control feedback and $x_k \in \mathbb{R}^n$ is the output used to evaluate the convergence rate.

For example, consider the following well-known first-order optimization algorithm:
\begin{equation}
    \begin{aligned}
    x_{k+1}=x_k+\gamma(x_k-x_{k-1})-\alpha \nabla f(y_k)\\
    y_k=x_k+\eta(x_k-x_{k-1})
    \end{aligned}
\end{equation}
where $\alpha, \gamma$ and $\eta$ are non-negative scalars. By defining the state vector $\xi_k=[x_{k-1}^{\top},  x_{k}^{\top}]^{\top}\in  \mathbb{R}^{2n}$, we can represent formula (16) in the form of (15), where the matrices ($A, B, C, E$) are given by
\begin{equation}
    \begin{aligned}
     &   A=\left[\begin{array}{cc}
0 & I_n \\ 
 -\gamma I_n & (1+\gamma)I_n
\end{array}\right]
\quad B=\left[\begin{array}{c}
0 \\
-\alpha I_n
\end{array}\right]\\
 & C=[-\eta I_n, (\eta+1)I_n]\quad \quad E=[0, I_n].
    \end{aligned}
\end{equation}
By choosing different values of $\gamma$ and $\eta$, we can formulate those typical GD methods, such as the conventional gradient decent method (where $\gamma=\eta=0$), the Heavy-ball method of Polyak\cite{b25} ($\eta=0$), and the Nesterov’s accelerated method ($\gamma=\eta \neq 0$). 

\section{Main Results}
We consider that $f: \mathbb{R}^n \rightarrow \mathbb{R} \cup\{+\infty\}$ is a $L-$smooth and $m-$strongly-convex function. 
Our objective is to solve the optimization problem:
\begin{equation}
    \min_{x\in R^n} f(x) ,
\end{equation}
and $x_*$ is the unique minimum of $f(x)$  and $\nabla f(x_*)=0$.

\textbf{Definition 4\cite{b26}:} A sequence $\{x_k\}$  is said to converge $R$-linearly to $x_*$
 with rate $\rho \in (0,1)$ if there is a constant $c>0$ such that
\begin{equation}
    ||x_k-x_*||_2\leq c\rho^k\quad \forall {k \in N_+}.
\end{equation}

\subsection{Adaptive fractional order gradient descent}
We propose the following adaptive fractional order gradient descent method (AFOGD):
    \begin{equation}
        {x}_{k+1}={x}_{k}-\alpha \nabla f\left({x}_{k}\right) \cdot \beta_k
        (\left\|x_k-x_{k-1}\right\|_2+\delta)^{1-\mu}
    \end{equation}
    and we require
    \begin{equation}
        0<c_1\leq \beta_k
        (\left\|x_k-x_{k-1}\right\|_2+\delta)^{1-\mu}\leq c_2 <+\infty ,
    \end{equation}
    where $0<\mu<2$, $c_1$ and $c_2$ are constants.
The formula (20) can be written as a dynamic system (15) with 
\begin{equation}
\begin{aligned}
     A=I_n\quad\quad B=-\alpha I_n \quad\quad C=I_n \quad\quad E=I_n\\
     u_k=\nabla f\left({x}_{k}\right) \cdot \beta_k
        (\left\|x_k-x_{k-1}\right\|_2+\delta)^{1-\mu}.
\end{aligned}
\end{equation}

\textbf{THEOREM 1:}
 If $f(x)$ is a $L-$smooth and $m-$strongly convex function, through the AFOGD method to solve the problem (18):\\
1. Suppose that the fixed points $\left(\xi_{\star}, u_{\star}, y_{\star}, x_{\star}\right)$  satisfy
\begin{equation}
    \xi_{\star}=A \xi_{\star}+B u_{\star},~y_{\star}=C \xi_{\star},~u_{\star}=\phi\left(y_{\star}\right), ~x_{\star}=E \xi_{\star}=y_{\star}
\end{equation}
for all $k$.\\
2. Suppose there exists symmetric matrix $ N$ such that the following inequalities hold for all $k$ :
\begin{equation}
    \begin{aligned}
0 & \leq e_k^{\top} N e_k,
\end{aligned}
\end{equation}
where $e_k=\left[\left(\xi_k-\xi_{\star}\right)^{\top}, \left(u_k-u_{\star}\right)^{\top}\right]^{\top}\in\mathbb{R}^{2n}$ .\\
3. Suppose there exists a non-negative constant $h$, a constant $\rho \in (0,1)$, a positive constant $\alpha$, and a positive definite matrix $P$ satisfying
\begin{equation}
    \quad M+h N \preceq 0 ,
\end{equation}
 where
 \begin{equation}
     M=\left[\begin{array}{cc}
A^{\top} P A-\rho^2 P & A^{\top} P B \\
B^{\top} P A & B^{\top} P B
\end{array}\right] .
 \end{equation}
Then the sequence $\left\{x_k\right\}$ $R-$linearly  converges to $x_*$.

\begin{proof}
We can select  a Lyapunov function
\begin{equation}
    V(\xi)=(\xi-\xi_*)^TP(\xi-\xi_*).
\end{equation}
Then,
    \begin{equation}
\begin{aligned}
   & V(\xi_{k+1})-\rho^2V(\xi_k)\\
   & =(\xi_{k+1}-\xi_*)^TP(\xi_{k+1}-\xi_*)-\rho^2(\xi_{k}-\xi_*)^TP(\xi_{k}-\xi_*)\\
   & =\left[\begin{array}{c}
\xi_k-\xi_* \\
u_k-u_*
\end{array}\right]^{\top}\left[\begin{array}{cc}
A^{\top} P A-\rho^2 P & A^{\top} P B \\
B^{\top} P A & B^{\top} P B
\end{array}\right] \left[\begin{array}{c}
\xi_k-\xi_* \\
u_k-u_*
\end{array}\right]\\
&= e_k^{\top} M e_k .
\end{aligned}
\end{equation}
From the above condition, we can know
\begin{equation}
     V(\xi_{k+1})-\rho^2V(\xi_k)\leq -h e_k^{\top} N e_k \leq 0 .
\end{equation}
Hence, 
\begin{equation}
     V(\xi_{k})=(\xi_{k}-\xi_*)^TP(\xi_{k}-\xi_*)\leq \rho^{2k}V(\xi_0) .
\end{equation}
Then $\lambda_0 ||\xi_{k}-\xi_{*}||_2^2\leq \rho^{2k}V(\xi_0)$, where $\lambda_0$ is the smallest eigenvalue of $P$.
This will lead to
\begin{equation}
\begin{aligned}
    & ||x_k-x_*||_2^2=||E(\xi_k-\xi_*)||_2^2\\
    & \leq ||E||^2\cdot||\xi_{k}-\xi_{*}||_2^2\leq \rho^{2k} \frac{V(\xi_0)||E|| }{\lambda_0} .
\end{aligned}
\end{equation}
Then take the square root of both sides and  get
\begin{equation}
    ||x_k-x_*||_2\leq \rho^k\sqrt{\frac{V(\xi_0)}{\lambda_0}}||E|| .
\end{equation}
\end{proof}

According to Theorem 1, if the AFOGD method proposed above is going to achieve  $R-$linear convergence, it needs to satisfy the above three conditions. For the first condition, it is obviously true  with $\xi_*=x_*=y_*$ and $u_*=0$. 

For the second requirement, we know that $f(x) $ is $m-$strong and $L-$smooth convex function so $f(x)$ satisfies inequality $(5)$. Since AFOGD has the adaptive step and meets  inequality $(21)$, then
\begin{equation}
    \begin{aligned}
    &\frac{m L}{m+L}\|x_k-x_*\|_2^2+\\
   & \frac{1}{(m+L)c_2}\|\nabla f(x_k)\cdot \beta_k
        (\left\|x_k-x_{k-1}\right\|_2+\delta)^{1-\mu} \|_2^2\\
       & \leq 
        \frac{1}{c_1}(\nabla f(x_k)\cdot \beta_k
        (\left\|x_k-x_{k-1}\right\|_2+\delta)^{1-\mu})^{\top}(x_k-x_*).
    \end{aligned}
\end{equation}
 This is equivalent to
 \begin{equation}
     \left[\begin{array}{c}
x_k-x_* \\
u_k-u_*
\end{array}\right]^{\top} N\left[\begin{array}{c}
x_k-x_* \\
u_k-u_*
\end{array}\right] \geq 0 ,
 \end{equation}
 where
 \begin{equation}
     N=\left[\begin{array}{cc}
\frac{-m L}{m+L} I_n &  \frac{1}{2c_1}I_n \\
 \frac{1}{2c_1}I_n & \frac{-1}{(m+L)c_2} I_n
\end{array}\right].
 \end{equation}
So the second condition can be achieved from the above analysis.

For the LMI $(25)$, we can select the $P=pI_n$ with $p>0$. Substitute $A, B$ and $  P$ into $(25)$:
\begin{equation}
    \left[\begin{array}{cc}
p_0-\rho^2p_0 & -\alpha p_0 \\
-\alpha p_0 & \alpha^2p_0
\end{array}\right] \otimes I_n +  h \left[\begin{array}{cc}
\frac{-m L}{m+L}  &  \frac{1}{2c_1} \\
\frac{1}{2c_1}  & \frac{-1}{(m+L)c_2} 
\end{array}\right] \otimes I_n \preceq 0
\end{equation}
and then we can reduce the dimension and get the  following ``$2\times$'' LMI:
\begin{equation}
    \left[\begin{array}{cc}
p_0-\rho^2p_0 & -\alpha p_0 \\
-\alpha p_0 & \alpha^2p_0
\end{array}\right]  +  h \left[\begin{array}{cc}
\frac{-m L}{m+L}  & \frac{1}{2c_1}  \\
\frac{1}{2c_1}  & \frac{-1}{(m+L)c_2} 
\end{array}\right] \preceq 0.
\end{equation}
Without loss of generality, we can select $p_0=1$, so the LMI becomes:
\begin{equation}
     \left[\begin{array}{cc}
1-\rho^2-\frac{h mL}{(m+L)} & -\alpha+\frac{h}{2c_1} \\
-\alpha+\frac{h}{2c_1} & \alpha^2-\frac{h}{(m+L)c_2}
\end{array}\right] \preceq 0.
\end{equation}
To simplify the calculation, we can select $h=2\alpha c_1$. In order for the LMI to be true, it is equivalent to:
\begin{equation}
    1-\rho^2-\frac{2\alpha c_1 mL}{m+L}\leq 0\quad and \quad \alpha^2-\frac{2\alpha c_1}{(m+L)c_2} \leq 0.
\end{equation}
Therefore, it can be found that when $0<\alpha \leq \frac{2c_1}{(m+L)c_2}$ and $ 1-\frac{2\alpha c_1 mL}{m+L} \leq \rho^2<1$, the LMI will be true.

In summary, through the above analysis, all three conditions of Theorem 1 can be satisfied. So AFOGD can ensure $R-$linear convergence for strongly convex functions.

As a summary of this subsection, the AFOGD method is formulated as the following Algorithm \ref{algorithm1}.

\begin{algorithm}
    \caption{AFOGD}\label{algorithm1}
    \SetKwInOut{Input}{Input}\SetKwInOut{Output}{Output}
    \Input{initial point $x_0$ and $x_1$, $\delta>0$, $\epsilon>0$, $0<\mu<2$, $c_1$, $c_2$ and let $0<\alpha \leq \frac{2c_1}{(m+L)c_2}$, $k=0$, $k_{max}$}
    \Output{optimal $x_*$ and $f(x_*)$}
    While
    {$k<k_{max}$ $ \&\& $ $\nabla f(x_k)\geq \epsilon$ do}\\
    Find $\beta_k$: $c_1\leq \beta_k
        (\left\|x_k-x_{k-1}\right\|_2+\delta)^{1-\mu}\leq c_2$\\
    Compute: ${x}_{k+1}={x}_{k}-\alpha \nabla f\left({x}_{k}\right) \cdot \beta_k
        (\left\|x_k-x_{k-1}\right\|_2+\delta)^{1-\mu}$\\
    Update: $k=k+1$\\
\end{algorithm}

\subsection{Adaptive fractional order accelerated gradient descent}
Subsequently, we propose the following adaptive fractional order accelerated gradient descent (AFOAGD)  algorithm:
\begin{equation}
    \begin{aligned}
    & {x}_{k+1}={y}_{k}-\alpha \nabla f\left({y}_{k}\right)\cdot\beta_k
        (\left\|y_k-y_{k-1}\right\|_2+\delta)^{1-\mu}\\
    & y_k=x_k+\eta(x_k-x_{k-1})
    \end{aligned}
\end{equation}
with an adaptive step:
\begin{equation}
        0<c_1\leq \beta_k(\left\|y_k-y_{k-1}\right\|_2+\delta)^{1-\mu}\leq c_2< +\infty .
\end{equation}
 By dynamical system (15), we select
 \begin{equation}
     \xi_k=[x_{k-1},x_{k}]^{\top},~u_k= \nabla f\left({y}_{k}\right)\cdot\beta_k
        (\left\|y_k-y_{k-1}\right\|_2+\delta)^{1-\mu}.
 \end{equation}
    Then we can get easily
\begin{equation}
    \begin{aligned}
    A=\left[\begin{array}{cc}
0 & I_n \\ 
 -\eta I_n & (1+\eta)I_n
\end{array}\right]
\quad B=\left[\begin{array}{c}
0 \\
-\alpha I_n
\end{array}\right]
    \end{aligned}
\end{equation}
and
\begin{equation}
    C=[-\eta I_n,(\eta+1)I_n]\quad \quad E=[0,I_n] .
\end{equation}

\textbf{Lemma 1:}  Let $f$ be $L-$smooth and $m-$strongly convex with
$m > 0$, and $x_*$ be the unique point satisfying $\nabla f(x_*)=0$. If the AFOAGD is applied, then the following inequality holds for all trajectories:
\begin{equation}
     f(x_{k+1})-f(x_*)\leq e_k^{\top}N^1e_k ,
\end{equation}
where
\begin{equation}
 N^1=\left[\begin{array}{ccc}
\frac{-\eta^2m}{2} I_n & \frac{\eta(\eta+1)m}{2} I_n & -\frac{\eta}{2c_1} I_n \\
\frac{\eta(\eta+1)m}{2} I_n & \frac{-(\eta+1)^2m}{2} I_n & \frac{\eta+1}{2c_1} I_n\\
-\frac{\eta}{2c_1} I_n & \frac{\eta+1}{2c_1} I_n & (\frac{1}{2}\alpha^2L-\frac{\alpha}{c_2})I_n
\end{array}\right] .
\end{equation}

\begin{proof}
Because $f(x)$ is $L-$smooth, then
\begin{equation}
    f(x_{k+1})-f(y_k)\leq \nabla f(y_k)^{\top}(x_{k+1}-y_k)+\frac{L}{2}||x_{k+1}-y_k||_2^2 .
\end{equation}
We notice that:
\begin{equation}
\begin{aligned}
 &\nabla f(y_k)^{\top}(x_{k+1}-y_k)\\& =-\alpha \nabla f(y_k)^{\top} \nabla f\left({y}_{k}\right)\cdot \beta_k(\left\|y_k-y_{k-1}\right\|_2+\delta)^{1-\mu} \leq 0 .
\end{aligned}
\end{equation}
Hence, according to $(41)$, we can get
\begin{equation}
\begin{aligned}
  & f(x_{k+1})-f(y_k)\\
  & \leq \frac{1}{c_2} \beta_k(\left\|y_k-y_{k-1}\right\|_2+\delta)^{1-\mu}\cdot \nabla f(y_k)^{\top}(x_{k+1}-y_k)\\
  & +\frac{L}{2}||x_{k+1}-y_k||_2^2 .
\end{aligned}
\end{equation}
The inequality $(49)$ can be rewritten as
\begin{equation}
\begin{aligned}
 &f(x_{k+1})-f(y_k)\leq\\
 & \left[\begin{array}{c}
x_{k+1}-y_k \\
u_k-u_*
\end{array}\right]^{\top}\left[\begin{array}{cc}
\frac{L}{2}I_n & \frac{1}{2c_2}I_n \\
\frac{1}{2c_2}I_n & 0
\end{array}\right] \left[\begin{array}{c}
 x_{k+1}-y_k\\
u_k-u_*
\end{array}\right]
\end{aligned}
\end{equation}
Based on the dynamical system (15), we have
\begin{equation}
    \left[\begin{array}{c}
 x_{k+1}-y_k\\
u_k-u_*
\end{array}\right]=\left[\begin{array}{cc}
EA-C & EB \\
0 & I_n
\end{array}\right]  \left[\begin{array}{c}
 \xi_{k}-\xi_*\\
u_k-u_*
\end{array}\right].
\end{equation}
Substituting (43), (44) and (51) into (50), we can write
\begin{equation}
    f(x_{k+1})-f(y_k)\leq e_k^{\top}W^1e_k ,
\end{equation}
where  
\begin{equation}
W^1=\left[\begin{array}{ccc}
0 & 0 & 0 \\
0 & 0 & 0\\
0 & 0 & (\frac{1}{2}\alpha^2L-\frac{\alpha}{c_2})I_n
\end{array}\right] .
\end{equation}
In addition, because $f$ is $m-$strongly convex, we know
    \begin{equation}
        f(y_k)-f(y_*)\leq \nabla f(y_k)^{\top}(y_k-y_*)-\frac{m}{2}||y_k-y_*||_2^2.
    \end{equation}
    We notice that 
    \begin{equation}
        \nabla f(y_k)^{\top}(y_k-y_*)\geq 0 ,
    \end{equation}
    then based on (41), we can get
    \begin{equation}
    \begin{aligned}
      & f(y_k)-f(y_*)\leq \\ & \frac{1}{c_1}\beta_k(\left\|y_k-y_{k-1}\right\|+\delta)^{1-\mu}\cdot\nabla f(y_k)^{\top}(y_{k}-y_*)\\&
      -\frac{m}{2}||y_{k}-y_*||_2^2 .
    \end{aligned}
    \end{equation}
 The inequality (56) can be rewritten as
\begin{equation}
\begin{aligned}
& f(y_{k})-f(y_*)\leq \\
   & \left[\begin{array}{c}
y_{k}-y_* \\
u_k-u_*
\end{array}\right]^{\top}\left[\begin{array}{cc}
-\frac{m}{2}I_n & \frac{1}{2c_1}I_n \\
\frac{1}{2c_1}I_n & 0
\end{array}\right] \left[\begin{array}{c}
 y_{k}-y_*\\
u_k-u_*
\end{array}\right].
\end{aligned}
\end{equation}
Once again, recalling the dynamical system (15), 
\begin{equation}
    \left[\begin{array}{c}
 y_{k}-y_*\\
u_k-u_*
\end{array}\right]=\left[\begin{array}{cc}
C & 0 \\
0 & I_n
\end{array}\right]  \left[\begin{array}{c}
 \xi_{k}-\xi_*\\
u_k-u_*
\end{array}\right].
\end{equation}
Substituting (43), (44) and (58) into (57), we can get
\begin{equation}
    f(y_{k})-f(y_*)\leq e_k^{\top}W^2e_k ,
\end{equation}
where
\begin{equation}
\begin{aligned}
  W^2=\left[\begin{array}{ccc}
-\frac{\eta^2m}{2}I_n & \frac{\eta(\eta+1)m}{2}I_n & -\frac{\eta}{2c_1} I_n \\
\frac{\eta(\eta+1)m}{2}I_n & -\frac{(\eta+1)^2m}{2}I_n & \frac{\eta+1}{2c_1}I_n\\
 -\frac{\eta}{2c_1} I_n& \frac{\eta+1}{2c_1}I_n & 0
\end{array}\right] 
\end{aligned}.
\end{equation}
By adding (52) to (59), we get
\begin{equation}
    f(x_{k+1})-f(x_*)\leq e_k^{\top}N^1e_k .
\end{equation}
\end{proof}

 \textbf{Lemma 2:} Let $f$ be $L-$smooth and $m-$strongly convex with
$m > 0$, and $x_*$ be the unique point satisfying $\nabla f(x_*)=0$. Considering the AFOAGD, the following inequalities holds for all trajectories :
\begin{equation}
    f(x_{k+1})-f(x_k)\leq e_k^{\top}N^2e_k ,
\end{equation}
where
\begin{equation}
\begin{aligned}
     & N^2=\psi(\nabla f(y_k)^{\top}(y_k-x_k)) \cdot\\
    & \left[\begin{array}{ccc}
\frac{-\eta^2m}{2} I_n & \frac{\eta^2m}{2} I_n & -\frac{\eta}{2c_1} I_n \\
\frac{\eta^2m}{2} I_n & \frac{-\eta^2m}{2} I_n & \frac{\eta}{2c_1} I_n\\
-\frac{\eta}{2c_1} I_n & \frac{\eta}{2c_1} I_n & (\frac{1}{2}\alpha^2L-\frac{\alpha}{c_2})I_n
\end{array}\right]\\
+& (1-\psi(\nabla f(y_k)^{\top}(y_k-x_k))\cdot\\
&\left[\begin{array}{ccc}
\frac{-\eta^2m}{2} I_n & \frac{\eta^2m}{2} I_n & -\frac{\eta}{2c_2} I_n \\
\frac{\eta^2m}{2} I_n & \frac{-\eta^2m}{2} I_n & \frac{\eta}{2c_2} I_n\\
-\frac{\eta}{2c_2} I_n & \frac{\eta}{2c_2} I_n & (\frac{1}{2}\alpha^2L-\frac{\alpha}{c_2})I_n
\end{array}\right],
\end{aligned}
\end{equation}
with $\psi $ is a function: $\mathbb{R}\rightarrow \mathbb{R}$ and satisfies that if $x\geq0, ~ \psi(x)=1$, otherwise $\psi(x)=0$.

\begin{proof}
Because $f$ is a $m-$strongly convex function,
\begin{equation}
        f(y_k)-f(x_k)\leq \nabla f(y_k)^{\top}(y_k-x_k)-\frac{m}{2}||y_k-x_k||_2^2.
    \end{equation}
    \textbf{Case 1:} $\nabla f(y_k)^{\top}(y_k-x_k)\geq 0$.\\
    In this case, from inequalities (41), we can get
    \begin{equation}
    \begin{aligned}
     & f(y_{k})-f(x_k)\leq\\ & \frac{1}{c_1}\beta_k(\left\|y_k-y_{k-1}\right\|_2+\delta)^{1-\mu}\cdot \nabla f(y_k)^{\top}(y_{k}-x_k)\\
     &-\frac{m}{2}||y_{k}-x_k||_2^2 .
       \end{aligned}
\end{equation}
The inequality (65) can be rewritten as
\begin{equation}
\begin{aligned}
 &f(y_{k})-f(x_k)\leq\\
   & \left[\begin{array}{c}
y_{k}-x_k \\
u_k-u_*
\end{array}\right]^{\top}\left[\begin{array}{cc}
-\frac{m}{2}I_n & \frac{1}{2c_1}I_n \\
\frac{1}{2c_1}I_n & 0
\end{array}\right] \left[\begin{array}{c}
 y_{k}-x_k\\
u_k-u_*
\end{array}\right].
\end{aligned}
\end{equation}
By dynamical system (15),
\begin{equation}
    \left[\begin{array}{c}
 y_{k}-x_k\\
u_k-u_*
\end{array}\right]=\left[\begin{array}{cc}
C-E & 0 \\
0 & I_n
\end{array}\right]  \left[\begin{array}{c}
 \xi_{k}-\xi_*\\
u_k-u_*
\end{array}\right].
\end{equation}
 Substituting (43), (44) and (67) into (66), we can get
\begin{equation}
     f(y_{k})-f(x_k)\leq e_k^{\top}W^3e_k ,
\end{equation}
where
\begin{equation}
\begin{aligned}
 W^3=\left[\begin{array}{ccc}
\frac{-\eta^2m}{2} I_n & \frac{\eta^2m}{2} I_n & -\frac{\eta}{2c_1} I_n \\
\frac{\eta^2m}{2} I_n & \frac{-\eta^2m}{2} I_n & \frac{\eta}{2c_1} I_n\\
-\frac{\eta}{2c_1} I_n & \frac{\eta}{2c_1} I_n & 0
\end{array}\right].
\end{aligned}
\end{equation}
\textbf{Case 2:} $\nabla f(y_k)^{\top}(y_k-x_k)< 0$.\\
    Here, based on (41), we can get
\begin{equation}
\begin{aligned}
     & f(y_{k})-f(x_k)\leq\\ & \frac{1}{c_2}\beta_k(\left\|y_k-y_{k-1}\right\|_2+\delta)^{1-\mu}\cdot 
     \nabla f(y_k)^{\top}(y_{k}-x_k)\\
     &-\frac{m}{2}||y_{k}-x_k||_2^2 .
\end{aligned}
\end{equation}
The above inequality can be rewritten as
\begin{equation}
\begin{aligned}
 &f(y_{k})-f(x_k)\leq\\
 & \left[\begin{array}{c}
y_{k}-x_k \\
u_k-u_*
\end{array}\right]^{\top}\left[\begin{array}{cc}
-\frac{m}{2}I_n & \frac{1}{2c_2}I_n \\
\frac{1}{2c_2}I_n & 0
\end{array}\right] \left[\begin{array}{c}
 y_{k}-x_k\\
u_k-u_*
\end{array}\right]
\end{aligned}
\end{equation}
 Substituting (43), (44) and into (71), we have
\begin{equation}
     f(y_{k})-f(x_k)\leq e_k^{\top}W^4e_k,
\end{equation}
where
\begin{equation}
\begin{aligned}
  W^4=\left[\begin{array}{ccc}
\frac{-\eta^2m}{2} I_n & \frac{\eta^2m}{2} I_n & -\frac{\eta}{2c_2} I_n \\
\frac{\eta^2m}{2} I_n & \frac{-\eta^2m}{2} I_n & \frac{\eta}{2c_2} I_n\\
-\frac{\eta}{2c_2} I_n & \frac{\eta}{2c_2} I_n & 0
\end{array}\right].
\end{aligned}
\end{equation}
By adding (52) to (68) and (72), we can get 
\begin{equation}
    f(x_{k+1})-f(x_k)\leq e_k^{\top}N^2e_k.
\end{equation}
\end{proof}

\noindent
From Lemma $1$ and Lemma $2$ above, we can obtain the following Theorem $2$.\\

\textbf{Theorem 2:}
$f(x)$ is a $L-$smooth and $m-$strongly convex function, and the AFOAGD method is applied to solve the problem (18). If there exists a $P$ ($P \succeq 0$), a non-negative constant $h$ and a constant $\rho \in (0,1)$ such that  the following LMI is true:
 \begin{equation}
      M+(1-\rho^2)N^1+\rho^2 N^2+hN^3 \preceq 0 ,
 \end{equation}
where
\begin{equation}
    M=\left[\begin{array}{cc}
A^{\top} P A-\rho^2 P & A^{\top} P B \\
B^{\top} P A & B^{\top} P B
\end{array}\right] 
\end{equation}
and 
\begin{equation}
    N^3=\left[\begin{array}{cc}
C & 0 \\
0 & I_n
\end{array}\right]^{\top} \left[\begin{array}{cc}
\frac{-m L}{m+L} I_n &  \frac{1}{2c_1}I_n\\
 \frac{1}{2c_1}I_n & \frac{-1}{(m+L)c_2} I_n
\end{array}\right] \left[\begin{array}{cc}
C & 0 \\
0 & I_n
\end{array}\right] ,
\end{equation}
then  the sequence $\left\{x_k\right\}$ $R-$linearly converges to $x_*$.
\begin{proof}
We select a Lyapunov function
\begin{equation}
    V(\xi)=(\xi-\xi_*)^TP(\xi-\xi_*)+f(x_k)-f(x_*).
\end{equation}
Then 
\begin{equation}
\begin{aligned}
   & V(\xi_{k+1})-\rho^2V(\xi_k)
   \\& =(\xi_{k+1}-\xi_*)^TP(\xi_{k+1}-\xi_*)-\rho^2(\xi_{k}-\xi_*)^TP(\xi_{k}-\xi_*)
   \\ & +f(x_{k+1})-f(x_*)-\rho^2(f(x_k)-f(x_*))
   \\ & =\left[\begin{array}{c}
\xi_k-\xi_* \\
u_k-u_*
\end{array}\right]^{\top}\left[\begin{array}{cc}
A^{\top} P A-\rho^2 P & A^{\top} P B \\
B^{\top} P A & B^{\top} P B
\end{array}\right] \left[\begin{array}{c}
\xi_k-\xi_* \\
u_k-u_*
\end{array}\right]\\& +\rho^2(f(x_{k+1})-f(x_k))+(1-\rho^2)(f(x_{k+1})-f(x_*))\\ &= e_k^{\top} M e_k+\rho^2(f(x_{k+1})-f(x_k))\\
& +(1-\rho^2)(f(x_{k+1})-f(x_*)).
\end{aligned}
\end{equation}
From Lemma 1 and Lemma 2 we can obtain
\begin{equation}
\begin{aligned}
     & V(\xi_{k+1})-\rho^2V(\xi_k)\\ &\leq e_k^{\top}(M+(1-\rho^2)N^1+\rho^2 N^2)e_k\leq -he_k^{\top}N^3e_k.
     \end{aligned}
\end{equation}
From the above inequality (34) we know
\begin{equation}
    \begin{aligned}
\left[\begin{array}{c}
y_k-y_* \\
u_k-u_*
\end{array}\right]^{\top} N\left[\begin{array}{c}
y_k-y_* \\
u_k-u_*
\end{array}\right]\geq 0 .
\end{aligned}
\end{equation}
By (43) and (44),
\begin{equation}
    \left[\begin{array}{c}
 y_{k}-y_*\\
u_k-u_*
\end{array}\right]=\left[\begin{array}{cc}
C & 0 \\
0 & I_n
\end{array}\right]  \left[\begin{array}{c}
 \xi_{k}-\xi_*\\
u_k-u_*
\end{array}\right].
\end{equation}
Hence, we can derive
\begin{equation}
    0\leq e_k^{\top} N^3 e_k .
\end{equation}
Furthermore, we can get
\begin{equation}
     V(\xi_{k+1})-\rho^2V(\xi_k)\leq 0,
\end{equation}
which means $V(\xi_k)\leq \rho^{2k}V(\xi_0)$. We also know  $(\xi_k-\xi_*)^TP(\xi_k-\xi_*)\geq 0$, so it is obvious that
\begin{equation}
    f(x_k)-f(x_*) \leq \rho^{2k}((\xi_0-\xi_*)^{\top}P(\xi_0-\xi_*)+f(x_0)-f(x_*)).
\end{equation}
Lastly, $f(x)$ is $m-$stongly convex then $\frac{m}{2}||x_k-x_*||_2^2\leq f(x_k)-f(x_*)$, so
\begin{equation}
    ||x_k-x_*||_2\leq \rho^k\sqrt{\frac{2V(\xi_0)}{m}} .
\end{equation}
\end{proof}
\textbf{Remark:} When $f(x)$ is only a convex function, the AFOAGD  can still guarantee  that $f(x)$ asymptotically converges to  $f(x_*)$ from formula (84).

In summary, the proposed AFOAGD method is illustrated in the following Algorithm~\ref{algorithm2}.

\begin{algorithm}
    \caption{AFOAGD}\label{algorithm2}
    \SetKwInOut{Input}{Input}\SetKwInOut{Output}{Output}
    \Input{initial point $x_0$, $x_1$ and $y_0$, $c_1$ and $c_2$,  $\delta>0$, $\epsilon>0$, $0<\mu<2$, $\alpha>0$, $\eta>0$,  $k=1$ and $k_{max}$.}
    \Output{optimal $x_*$ and $f(x_*)$}
    While
    {$k<k_{max}$ $ \&\& $ $\nabla f(y_k)\geq \epsilon$ do}\\
    Find $\beta_k$: $ c_1\leq \beta_k
        (\left\|y_k-x_{y-1}\right\|_2+\delta)^{1-\mu}\leq c_2 $\\
    Compute: 
    $y_k=x_k+\eta(x_k-x_{k-1})$\\
    Compute:
    ${x}_{k+1}={y}_{k}-\alpha \nabla f\left({y}_{k}\right) \cdot \beta_k
        (\left\|y_k-y_{k-1}\right\|_2+\delta)^{1-\mu}$\\
        
    Update: $k=k+1$\\
\end{algorithm}

\section{Simulation}
In order to verify the effectiveness of our proposed methods - AFOGD and AFOAGD, we perform the following numerical simulations. 

\subsection{Simulation 1}
We still choose the quadratic function $f(x_1,x_2)=2x_1^2+3x_2^2+3$ from \textbf{Example 1} with $c_1=0.8$ , $c_2=1.3$ and the rest of the initial conditions remain the same. The simulation results are shown in Fig. 2. As can be seen from Fig. 2, FOGD does not converge to the minimum point, and compared with GD method, the AFOGD proposed by us has faster convergence speed and smaller error. This shows the effectiveness and superiority of AFOGD.
\begin{figure}
    \centering
    \includegraphics[width=\linewidth]{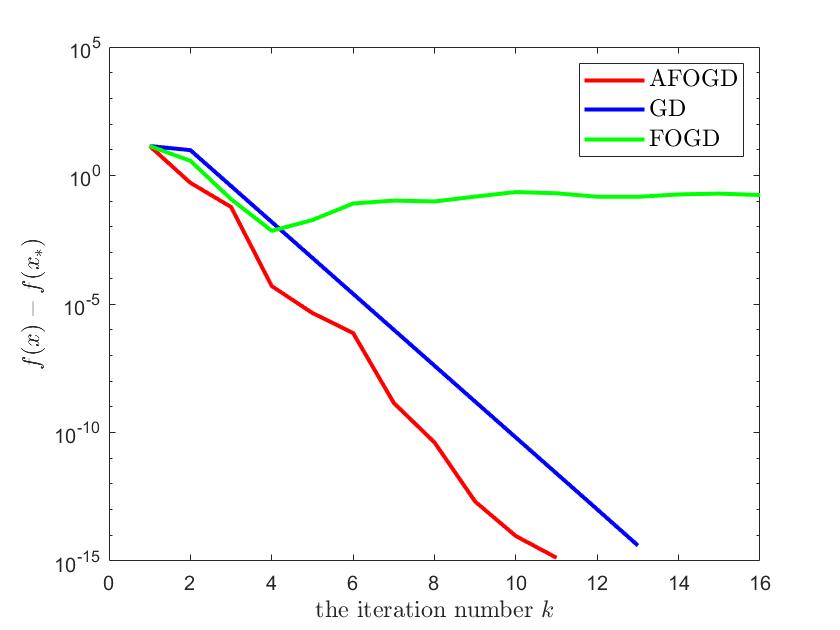}
    \caption{A numerial comparsion among AFOGD, GD and FOGD. The objective function is $f(x_1,x_2)=2x_1^2+3x_2^2+3$}
    \label{fig:my_label}
\end{figure}
\subsection{Simulation 2}
Consider a quadratic function $f(x_1,x_2)=8x_1^2+2x_2^2+4x_1+2x_2-1$, which is a smooth and strongly convex function with $L=8$ and $m=2$.
We use AFOAGD to solve the problem with $\alpha=0.1$, $\eta=0.2$, $c_1=0.5$ and $c_2=1$. Firstly, we need to solve the LMI (75) 
\begin{equation}
      M+(1-\rho^2)N^1+\rho^2 N^2+hN^3 \preceq 0.
 \end{equation}
When $\nabla f(y_k)^{\top}(y_k-x_k)\leq 0$, we can find  $h=0.2$ , $\rho^2=0.4$ and 
\begin{equation}
     \begin{aligned}
     P_2=\left[\begin{array}{cc}
4.1074 & -4.1697 \\
-4.1697 & 4.6191
\end{array}\right]
     \end{aligned}
\end{equation}
to make the above LMI  true. When $\nabla f(y_k)^{\top}(y_k-x_k)< 0$, $h=0.2$, $\rho^2=0.8$ and
\begin{equation}
     \begin{aligned}
     P_1=\left[\begin{array}{cc}
1.8345 & -1.6390 \\
-1.6390 & 7.0917
\end{array}\right]
     \end{aligned}
\end{equation}
can meet the LMI. So the LMI can be achieved. For $f(x_1,x_2)$, we can know the minimum point is $(-0.25, -0.5)^{\top}$ and the minimum is $-2$ easily. We select the initial point $x_0=[1.2,1.2]^{\top}$ and $y_0=x_1=(-1.12, 0.52)^{\top}$, the simulation result is showed Fig. 3. We compare AFOAGD with Heavy-ball ,GD and FOGD. As can be seen from the Fig. 3, FOGD does not converge to the minimum value, and compared with GD and Heavy-ball method, AFOAGD has smaller oscillation and faster convergence rate.
\begin{figure}
    \centering
    \includegraphics[width=\linewidth]{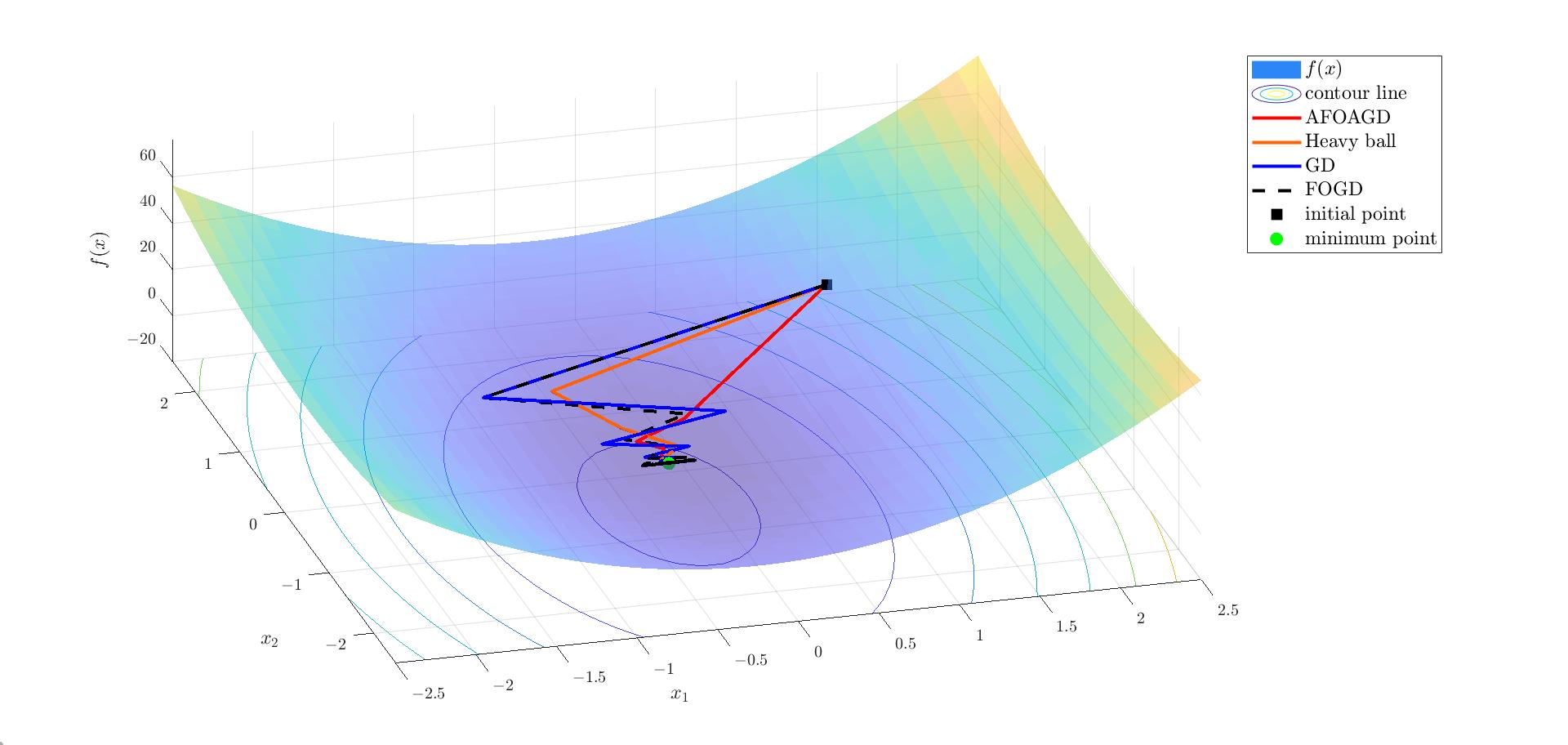}
    \caption{A numerical comparsion among AFOAGD, Heavy-ball method, GD and FOGD. The objective function is $f(x_1,x_2)=8x_1^2+2x_2^2+4x_1+2x_2-1$}
    \label{fig:my_label}
\end{figure}

\subsection{Simulation 3}
Consider 40 points in Fig. 4, the math model that we know in
advance is given by
\begin{equation}
    y=\theta_0+\theta_1x
\end{equation}
where $\theta_0$ and $\theta_1$ are parameters. We want to find the parameters that suit the points well. The error function can be defined as follow :
\begin{equation}
    J(\theta)=\frac{1}{80}\sum_{i=1}^{40}(\theta_0+\theta_1x_i-y_i)^2.
\end{equation}
Notice that $J(\theta)$ is only convex function so $m=0$, but the AFOAGD can still ensure that $J(\theta)$ achieves to reach minimum. In addition, we compare the NES algorithm with AFOAGD. We denote $\Theta=[\theta_0, \theta_1]^{\top}$ and we set the initial point $\Theta_0=[0.1, 0.8]^{\top}$ and $\Theta_1=[0, 0]^{\top}$, $\alpha=0.2$ and $\eta=0.1$, $c_1=1.3$ and $c_2=2$. By AFOAGD, we find the optimal parameter $\theta_0=0.4619$ and $\theta_1=1.9971$. And by NES, we find the optimal parameter $\theta_0=0.4617$ and $\theta_1=1.9970$. The fitting effect is shown in Fig. 4. In Fig. 5, we depict the changing process of $J(\theta)$.

\begin{figure}
    \centering
    \includegraphics[width=\linewidth]{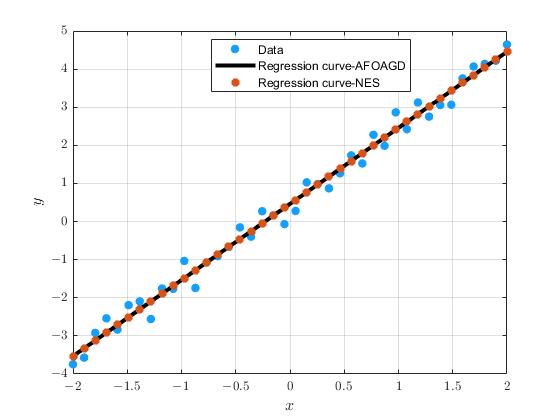}
    \caption{Linear regression using AFOAGD and NES}
    \label{fig:my_label}
\end{figure}

\begin{figure}
    \centering
    \includegraphics[width=\linewidth]{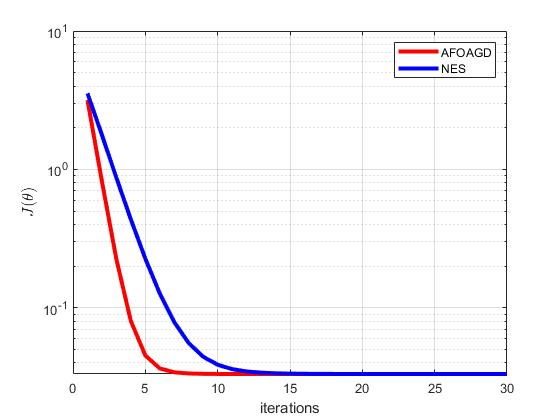}
    \caption{Curve of error}
    \label{fig:my_label}
\end{figure}

\section{CONCLUSIONS}

In this paper, inspired fractional order gradient descent may  diverge or converge to a region containing minimum point instead of converging to the true minimum point when the objective function is strongly convex function, two novel fractional order gradient descent methods  named AFOGD and  AFOAGD method are proposed.   Based on quadratic constraints and Lyapunov stability analysis from robust control theory, we certify that AFOGD and AFOAGD can achieve $R-$linear convergence when the objective function is $L-$smooth and $m-$strongly convex under some certain assumptions. Simulation 1 shows that AFOGD  can converge to the true minimum point for strongly convex function and compared with GD method, the AFOGD  has faster convergence speed and smaller error. In Simulation 2, we compare AFOAGD with Heavy-ball method ,GD and FOGD and show that AFOAGD has smaller oscillation and faster convergence rate. 
Moreover, we apply AFOAGD to linear regression and compare AFOAGD with Nesterov's accelerated method  in Simulation 3. Under the same initial conditions and parameter selection, the error curve obtained by AFOAGD algorithm approaches to zero more quickly. This shows that AFOGD and AFOAGD can be used in the fields of machine learning and system identification. Other theoretical analysis such the sensitivity of parameters or others in AFOGD and AFOAGD and their applications to the parameter identifications for some complex systems will be studied in the future.

\addtolength{\textheight}{-12cm}   




\end{document}